\documentclass[11pt]{article}
\usepackage{amsfonts,amsmath,amssymb,accents}
\usepackage{hyperref}
\usepackage{graphicx}

\newcommand{\status}{} 
\newcommand{\detail}[1]{\par\noi{\bf [Proof detail\ }{#1}
\hfill{\bf ]}\par\noi\hspace{-4pt}}
\renewcommand{\detail}[1]{}

\newcommand{\file}{wisk/rencon/birth.tex\quad}
\renewcommand{\file}{}

\newcommand{\dis}{\displaystyle}

\newcommand{\noi}{\noindent}

\newcommand{\halmos}{\rule{1ex}{1.4ex}}
\def \qed {\nopagebreak{\hspace*{\fill}$\halmos$\medskip}}
\newcommand{\med}{\medskip}

\newtheorem{theorem}{Theorem}[section]
\newtheorem{proposition}[theorem]{Proposition}
\newtheorem{corollary}[theorem]{Corollary}
\newtheorem{conjecture}[theorem]{Conjecture}
\newtheorem{lemma}[theorem]{Lemma}

\newtheorem{counterexample}[theorem]{Counterexample}
\newtheorem{remark}[theorem]{Remark}

\newcommand{\bt}{\begin{theorem}}
\newcommand{\et}{\end{theorem}}
\newcommand{\bl}{\begin{lemma}}
\newcommand{\el}{\end{lemma}}
\newcommand{\bp}{\begin{proposition}}
\newcommand{\ep}{\end{proposition}}
\newcommand{\bcor}{\begin{corollary}}
\newcommand{\ecor}{\end{corollary}}
\newcommand{\br}{\begin{remark}\rm}
\newcommand{\er}{\end{remark}}
\newcommand{\bcon}{\begin{conjecture}}
\newcommand{\econ}{\end{conjecture}}
\newcommand{\bcount}{\begin{counterexample}}
\newcommand{\ecount}{\end{counterexample}}

\newcommand{\be}{\begin{equation}}
\newcommand{\ee}{\end{equation}}
\newcommand{\ba}{\begin{array}}
\newcommand{\ea}{\end{array}}
\newcommand{\bc}{\be\begin{array}{r@{\,}c@{\,}l}}
\newcommand{\ec}{\end{array}\ee}

\newcommand{\de}{\delta}

\newcommand{\la}{\lambda}

\newcommand{\sig}{\sigma}

\newcommand{\Ei}{{\cal E}}
\newcommand{\Fi}{{\cal F}}

\newcommand{\R}{{\mathbb R}}
\newcommand{\N}{{\mathbb N}}

\renewcommand{\P}{{\mathbb P}}

\newcommand{\li}{\langle}
\newcommand{\re}{\rangle}

\newcommand{\dgg}{\dagger}

\begin{document}

%numbering formulas within sections
\makeatletter\@addtoreset{equation}{section}
\makeatother\def\theequation{\thesection.\arabic{equation}} 

%alternative layout for enumerate lists.
\renewcommand{\labelenumi}{{(\roman{enumi})}}

\title{\vspace{-3cm}Intertwining of birth-and-death processes}
\author{Jan~M.~Swart}
\date{{\small\file} \today}

\maketitle\vspace{-.8cm}
\status

\begin{abstract}\noi
It has been known for a long time that for birth-and-death processes started in zero the first passage time of a given level is distributed as a sum of independent exponentially distributed random variables, the parameters of which are the negatives of the eigenvalues of the stopped process. Recently, Diaconis and Miclo have given a probabilistic proof of this fact by constructing a coupling between a general birth-and-death process and a process whose birth rates are the negatives of the eigenvalues, ordered from high to low, and whose death rates are zero, in such a way that the latter process is always ahead of the former, and both arrive at the same time at the given level. In this note, we extend their methods by constructing a third process, whose birth rates are the negatives of the eigenvalues ordered from low to high and whose death rates are zero, which always lags behind the original process and also arrives at the same time.
\end{abstract}

\noi
%60J27	Continuous-time Markov processes on discrete state spaces
%
%15A18  Eigenvalues, singular values, and eigenvectors
%37A30  Ergodic theorems, spectral theory, Markov operators
%60G40  Stopping times; optimal stopping problems; gambling theory
%60J35  Transition functions, generators and resolvents
%60J80  Branching processes (Galton-Watson, birth-and-death, etc.)
{\it Acknowledgement.} Work sponsored by GA\v CR grant 201/09/1931.
%\vspace{12pt}

{\small\setlength{\parskip}{-2pt}\tableofcontents}
%\newpage

\section{Introduction}

\subsection{First passage times of birth-and-death processes}\label{S:tauN}

%\addtocounter{theorem}{-1}

Let $X=(X_t)_{t\geq 0}$ be the continuous-time Markov process in $\N=\{0,1,\ldots\}$, started from $X_0=0$, that jumps from $x-1$ to $x$ with {\em birth rate} $b_x>0$ and from $x$ to $x-1$ with {\em death rate} $d_x>0$ $(x\geq 1)$. Let
\be
\tau_N:=\inf\{t\geq 0:X_t=N\}\qquad(N\geq 1)
\ee
denote the first passage time of $N$. The following result has been known at least since \cite[Prop.~1]{KM59}.

\bp{\bf(Law of first passage times)}\label{P:tauN}
The first passage time $\tau_N$ is distributed as a sum of independent exponentially distributed random variables whose parameters $\la_1<\cdots<\la_N$ are the negatives of the nonzero eigenvalues of the generator of the process stopped in $N$.
\ep

Older proofs of this fact are based on a calculation of the Laplace transform of $\tau_N$ by purely algebraic methods, see \cite{DM09} for a historical overview. In the latter paper, Diaconis and Miclo gave for the first time a probabilistic proof of Proposition~\ref{P:tauN}, by coupling the process $X$ to another birth-and-death process $X^+$ with birth rates $b^+_1=\la_N,\ldots,b^+_N=\la_1$ and zero death rates, in such a way that $X_{t\wedge\tau_N}\leq X^+_t$ for all $t\geq 0$ and $X$ and $X^+$ arrive in $N$ at the same time. In the present paper, we will extend their methods by showing that $X$ and $X^+$ can in addition be coupled to a process $X^-$ with birth rates $b^+_1=\la_1,\ldots,b^+_N=\la_N$ and zero death rates, in such a way that $X^-_t\leq X_{t\wedge\tau_N}\leq X^+_t$ for all $t\geq 0$ and all three processes arrive in $N$ at the same time.

\subsection{Intertwining of Markov processes}\label{S:twine}

The coupling technique used by Diaconis and Miclo in \cite{DM09} is of a special kind, which is sometimes called {\em intertwining of Markov processes}. Let $X$ and $X'$ be continuous-time Markov processes with finite state spaces $S$ and $S'$ and generators $G$ and $G'$, respectively, and let $K$ be a probability kernel from $S$ to $S'$. Then $K$ defines a linear operator from $\R^{S'}$ to $\R^S$, also denoted by $K$, by the formula
\be
Kf(x):=\sum_{y\in S'}K(x,y)f(y).
\ee
The following result, which is based on an observation by Rogers and Pitman \cite{RP81}, was proved by Fill in \cite[Thm.~2]{Fil92}. (An independent proof can be found in \cite[Prop.~4]{AS10}).

\bp{\bf(Intertwining of Markov processes)}\label{P:twine} 
Assume that
\be\label{twine}
GK=KG'.
\ee
Then there exists a generator $\hat G$ of an $S\times S'$-valued Markov process with the property that if $(X,X')$ evolves according to $\hat G$ and satisfies
\be
\P[X'_0=y\,|\,X_0]=K(X_0,y)\qquad(y\in S'),
\ee
then
\be
\P[X'_t=y\,|\,(X_s)_{0\leq s\leq t}]=K(X_t,y)\qquad(t\geq 0,\ y\in S'),
\ee
and the processes $X$ and $X'$, on their own, are Markov processes evolving according to the generators $G$ and $G'$, respectively.
\ep
Algebraic relations of the type (\ref{twine}) are called intertwining relations, hence the name {\em intertwining of Markov processes}. We note that the operator $K$ needs in general not have an inverse, and even if it does, this inverse will in general not be associated to a probability kernel from $S'$ to $S$. In view of this, an intertwining of Markov processes is not a symmetric relation. To express this, following terminology introduced in \cite{AS10}, we will also say that in the set-up of Proposition~\ref{P:twine}, $X'$ is an {\em averaged Markov process on} $X$.

\subsection{Intertwining of birth-and-death processes}\label{S:main}

We are now ready to formulate our main result. Deviating slightly from our notation in Section~\ref{S:tauN}, we let $X=(X_t)_{t\geq 0}$ be a continuous-time Markov process with state space $\{0,\ldots,N\}$, started from $X_0=0$, that jumps from $x-1$ to $x$ with birth rate $b_x$ and from $x$ to $x-1$ with death rate $d_x$, where $b_1,\ldots,b_N>0$, $d_1,\ldots,d_{N-1}>0$, but $d_N=0$, i.e., $X$ is the stopped process from Section~\ref{S:tauN}. We let $G$ denote the generator of $X$, i.e.,
\be\label{Gdef}
Gf(x):=b_{x+1}\big(f(x+1)-f(x)\big)+d_x\big(f(x-1)-f(x)\big)
\qquad(0\leq x\leq N),
\ee
where $f:\{0,\ldots,N\}\to\R$ is a real function and we adopt the convention that $d_0=0$ and $b_{N+1}=0$ so that the corresponding terms in (\ref{Gdef}) are zero, regardless of the (fictive) values of $f$ in $-1$ and $N+1$. The following theorem is our main result.

\bt{\bf(Intertwining of birth-and-death processes)}\label{T:main}
The operator $G$ has $N+1$ distinct eigenvalues $0=-\la_0>-\la_1>\cdots>-\la_N$. Let $X^-$ and $X^+$ be the pure birth processes in $\{0,\ldots,N\}$, started from $X^-_0=X^+_0=0$, with birth rates $b^-_1:=\la_1,\ldots,b^-_N:=\la_N$ and $b^+_1:=\la_N,\ldots,b^+_N:=\la_1$, respectively, and let $G^-$ and $G^+$ be their generators. Then there exist probability kernels $K^-$ and $K^+$ on $\{0,\ldots,N\}$ satisfying
\be\ba{ll}\label{Kprop}
\dis K^-(x,\{0,\ldots,x\})=1,\quad&\dis K^+(x,\{0,\ldots,x\})=1,
\quad(0\leq x\leq N)\\[5pt]
\dis K^-(N,N)=1,\quad&\dis K^+(N,N)=1,
\ec
and
\be\label{Gpm}
{\rm(i)}\ \ K^+G=G^+K^+
\quad\mbox{and}\quad
{\rm(ii)}\ \ GK^-=K^-G^-.
\ee
Moreover, the processes $X^-,X$, and $X^+$ can be coupled in such a way that
\be\ba{rr@{\,}c@{\,}ll}\label{Xcoup}
{\rm(i)}&\dis\P[X_t=y\,|\,(X^+_s)_{0\leq s\leq t}]&=&\dis K^+(X^+_t,y)
\qquad&(t\geq 0,\ 0\leq y\leq N),\\[5pt]
{\rm(ii)}&\dis\P[X^-_t=y\,|\,(X^+_s,X_s)_{0\leq s\leq t}]&=&\dis K^-(X_t,y)
\qquad&(t\geq 0,\ 0\leq y\leq N).
\ec
\et
The existence of a kernel $K^+$ such that (\ref{Gpm})~(i) and (\ref{Xcoup})~(i) hold has been proved before in \cite[Prop.~10]{DM09}. Our new contribution is the construction of the kernel $K^-$ such that moreover (\ref{Gpm})~(ii) and (\ref{Xcoup})~(ii) hold. It is easy to see that formulas (\ref{Kprop}) and (\ref{Xcoup}) imply that
\be\ba{rl}
{\rm(i)}&\dis X^-_t\leq X_t\leq X^+_t\qquad(t\geq 0),\\[5pt]
{\rm(ii)}&\dis\tau^-_N=\tau_N=\tau^+_N,
\ec
where $\tau_N:=\inf\{t\geq 0:X_t=N\}$ and $\tau^-_N$ and $\tau^+_N$ are defined similarly for $X^-$ and $X^+$, respectively. We note that $X^-$ and $X^+$ move, in a sense, in the slowest resp.\ fastest possible way from $0$ to $N$, given that they have to arrive at exactly the same time as $X$. Note that, using terminology introduced at the end of Section~\ref{S:twine}, $X$ is an averaged Markov process on $X^+$ and $X^-$ is an averaged Markov process on $X$.

\subsection{Discussion}

In comparison to the paper by Diaconis and Miclo \cite{DM09}, the present paper does not add too much that is new. In particular the construction of the kernel $K^-$ in Theorem~\ref{T:main} is very similar to the construction of the kernel $K^+$, which was already carried out in \cite{DM09}. However, we believe that the observation that both constructions are possible, with an interesting symmetry between them, is of some interest.

The (new) construction with the process $X^-$ has in fact one advantage over the construction with $X^+$, since Proposition~\ref{P:twine} and formula (\ref{Gpm}) imply that the process $X$ started in {\em any} initial state can be coupled to a process $X^-$ with the same dynamics as in Theorem~\ref{T:main}, in such a way that $\P[X^-_t=y\,|\,(X_s)_{0\leq s\leq t}]=K^-(X_t,y)$ for all $0\leq y\leq N$ and $t\geq 0$. This implies that for a general initial state $X_0=x\in\{0,\ldots,N\}$, the stopping time $\tau_N$ is distributed as $\sum_{y=Z}^N\sig_y$ where $\sig_1,\ldots,\sig_N$ are independent exponentially distributed random variables with parameters $\la_1,\ldots,\la_N$ and $Z$ is an independent $\{0,\ldots,N\}$-valued random variable with law $K^-(x,\,\cdot\,)$. Note that the (old) coupling with the process $X^+$ forces one to start the process $X$ in an initial law that is a convex combination of the laws $K^+(x,\,\cdot\,)$ with $0\leq x\leq N$, hence no conclusions can be drawn for arbitrary initial states.

On the other hand, the methods of \cite{DM09} can also be used to study birth-and-death processes on $\{0,\ldots,N\}$ whose death rate $d_N$ is not zero and which, therefore, converge in law to a unique equilibrium. In particular, Diaconis and Miclo use a generalization of their intertwining relation (\ref{Gpm})~(i) to construct a fastest strong stationary time for such processes (we refer to \cite{DM09} for the definition). In contrast, it seems that the interwining relation (\ref{Gpm})~(ii) does not generalize to such a setting.

On a more general level, one may ask what the advantage is of a `probabilistic' proof of Proposition~\ref{P:tauN} as opposed to older, more algebraic proofs. Since most of the work behind Theorem~\ref{T:main} goes into proving the intertwining relations (\ref{Gpm}), one might even argue that the present proof is still rather algebraic in nature, although with a strong probabilistic flavour. In this context, it is interesting to note that the fact that $G$ is diagonalizable with real, distinct eigenvalues follows as a result of our proofs (in particular, this follows from a repeated application of the Perron-Frobenius theorem) and does not have to be provided by some extra argument (based on, for example, reversibility).

In general, diagonalizing a generator of a Markov process gives very strong information about the process, but in practice, if the state space is large, it is hard to get good information about the position of eigenvalues etc. The idea of interwining generators with transition kernels may in some cases be a good way to transform generators of complicated processes into generators of more simple processes and thus provide a more probabilistic alternative to diagonalization.

The methods of this paper can certainly be extended to one-dimensional processes with two traps, to dicrete-time processes, and to one-dimensional diffusions. Miclo \cite{Mic10} has proved a generalization of Proposition~\ref{P:tauN} for reversible Markov chains. In \cite{AS10}, intertwining relations were used to estimate the time to extinction for large hierarchical contact processes. The present work was partly motivated by an open problem from that paper. (To be precise, Question~$1^\circ$ from Section~3.3.)

%Method a priori not restricted to reversible processes but all depends on whether one can find a clever kernel.

\section{Proofs}

\subsection{Leading eigenvectors}

Let $X$ be the birth-and-death process in $S_N:=\{0,\ldots,N\}$ from Section~\ref{S:main} and let $G:\R^{S_N}\to\R^{S_N}$ be its generator, defined in (\ref{Gdef}). We equip $\R^{S_N}$ with the usual inner product $\li\pi|f\re:=\sum_{x=0}^N\pi(x)f(x)$ and let $G^\dgg$ denote the adjoint of $G$ with respect to this inner product. Then
\be\label{Gdagger}
G^\dgg\pi(x)=b_x\pi(x-1)-b_{x+1}\pi(x)+d_{x+1}\pi(x+1)-d_x\pi(x),
\ee
where as in (\ref{Gdef}) we use the convention that $d_0=0$ and $b_{N+1}=0$ so that the corresponding terms in (\ref{Gdagger}) are zero, regardless of the (fictive) values of $\pi$ in $-1$ and $N+1$.

Since $\de_N$ (the delta mass in $N$) is the unique invariant law of $X$, the eigenvalue $0$ of the generator $G$ has multiplicity one and its unique left and right eigenvectors are $\de_N$ and the constant function $1$, respectively. We will need the following result on the next largest eigenvalue and its left and right eigenvectors.

\bl{\bf(Leading eigenvectors)}\label{L:eigen}
There exists a $\la>0$ and $f,\pi\in\R^{S_N}$ such that
\[\ba{rl}
{\rm(i)}&\mbox{$f$ is strictly decreasing on $\{0,\ldots,N\}$ and satisfies $f(0)=1$, $f(N)=0$,}\\[5pt]
{\rm(ii)}&\mbox{$\pi$ is strictly positive on $\{0,\ldots,N-1\}$ and satisfies $\sum_{x=0}^{N-1}\pi(x)=1=-\pi(N)$,}\\[5pt]
{\rm(iii)}&Gf=-\la f\quad\mbox{and}\quad G^\dgg\pi=-\la\pi.
\ea\]
\el
{\bf Proof} Set
\be\ba{lll}
e(x):=\de_x\quad&(0\leq x\leq N-1)\quad&\mbox{and}\quad e(N):=1,\\
\xi(x):=\de_x-\de_N\quad&(0\leq x\leq N-1)\quad&\mbox{and}\quad \xi(N):=\de_N.
\ec
Then $\{e(0),\ldots,e(N)\}$ is a basis for $\R^{S_N}$ and $\{\xi(0),\ldots,\xi(N)\}$ is its associated dual basis, i.e., $\li e(x)|\xi(y)\re=1_{\{x=y\}}$. Set
\bc
\Ei&:=&{\rm span}\{e(0),\ldots,e(N-1)\}=\{f\in\R^{S_N}:f(N)=0\},\\[5pt]
\Fi&:=&{\rm span}\{\xi(0),\ldots,\xi(N-1)\}=\{\pi\in\R^{S_N}:\sum_{x=0}^N\pi(x)=0\}.
\ec
Since $N$ is a trap for the process $X$, it is easy to see that the operator $G$ maps the space $\Ei$ into itself. Since the coordinates of a vector in $\Ei$ with respect to the basis $\{e(0),\ldots,e(N)\}$ are the same as its coordinates with respect to the standard basis $\{\de_0,\ldots,\de_N\}$, it follows that with respect to the basis $\{e(0),\ldots,e(N)\}$, the matrix $[G]$ of $G$ has the form
\be
[G]=\left(\ba{cc}A&0\\ 0&0\ea\right),
\ee
where $A(x,y)=G(x,y)$ for $0\leq x,y\leq N-1$. The restriction of the process $X$ to the space $\{0,\ldots,N-1\}$ is irreducible in the sense that there is a positive probability of going from any state to any other state. Therefore, by applying the Perron-Frobenius theorem (see, e.g.\ Chapter XIII, \S 2, Theorem~2 in \cite{Gan00}) to $A+cI$ and its adjoint for some sufficiently large $c$, one finds that $A$ has a real eigenvalue $-\la$ of multiplicity one, which is larger than all other real eigenvalues, and associated left and right eigenvectors $\pi\in\Fi$ and $f\in\Ei$ that are strictly positive with respect to the bases $\{\xi(0),\ldots,\xi(N-1)\}$ and $\{e(0),\ldots,e(N-1)\}$, respectively. Since Markov semigroups are contractive we have $-\la\leq 0$ and since the eigenvalue zero of $G$ has multiplicity one and belongs to different left and right eigenvectors, we conclude that $-\la<0$. Since we can always normalize our eigenvectors such that $\sum_{x=0}^{N-1}\pi(x)=1$ and $\max_{x=0}^{N-1}f(x)=1$, this proves all statements of the lemma except for the fact that $f$ is strictly decreasing.

To prove this latter fact, we observe that by the facts that $Gf=-\la f$ and $f>0$ on $\{0,\ldots,N-1\}$,
\be
b_1\big(f(1)-f(0)\big)=-\la f(0)<0,
\ee
which show that $f(0)>f(1)$. By the same argument,
\be
b_{x+1}\big(f(x+1)-f(x)\big)=-\la f(x)-d_x\big(f(x-1)-f(x)\big)<0\qquad(1\leq x\leq N-1),
\ee
from which we see by induction that $f(x)>f(x+1)$ for all $0\leq x\leq N-1$.\qed

\subsection{Intertwining the fast process}

In this section, we prove the existence of a kernel $K^+$ satisfying (\ref{Kprop}) and (\ref{Gpm}). Our proof is basically the same as the proof given in \cite{DM09}, but as a preparation for the next section it will be convenient to review their proof and shorten it somewhat. The proof in \cite{DM09} is written in such a way as to make clear how the authors arrived at their argument and uses discrete derivatives that are presumably also useful if one wants to generalize the theory to one-dimensional diffusions. If our only aim is Theorem~\ref{T:main}, however, we can summarize their arguments quite a bit.

\begin{figure}
\begin{center}
\includegraphics[width=15cm]{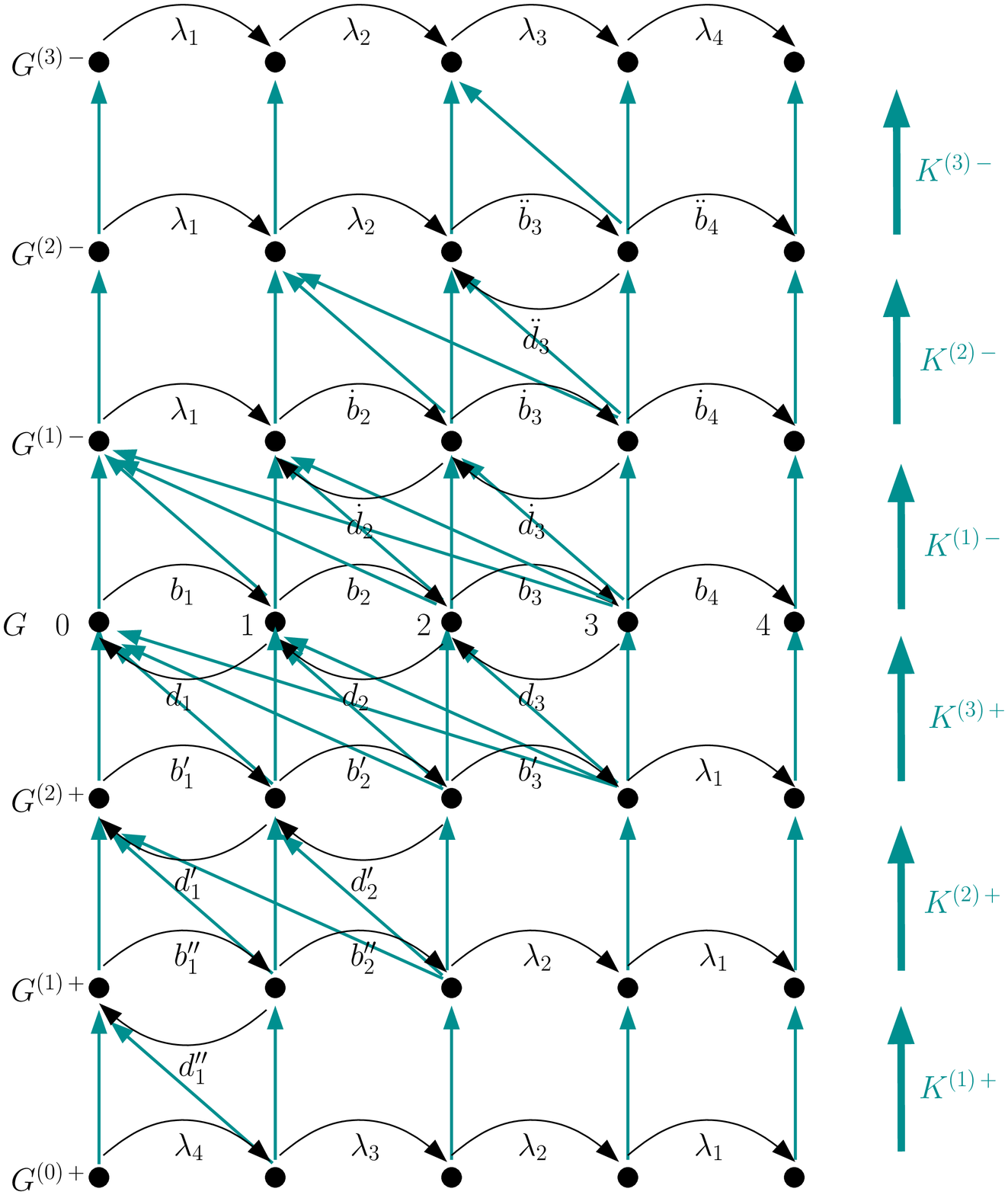}
\caption{Intertwining of birth and death processes. In this picture $N=4$. All nonzero transition rates and probabilities have been indicated with arrows.}
\label{fig:twine}
\end{center}
\end{figure}

The kernel $K^+$ will be constructed as the concatenation of an inductively defined sequence of kernels $K^{(N-1)\,+},\ldots,K^{(1)\,+}$. Associated with these kernels is a sequence of generators $G^{(N-1),+},\ldots,G^{(0)\,+}$ of birth-and-death processes in $\{0,\ldots,N\}$ satisfying the intertwining relations
\be\label{KGMplus}
K^{(M)\,+}G^{(M)\,+}=G^{(M-1)\,+}K^{(M)\,+}\qquad(1\leq M\leq N-1),
\ee
where the process with generator $G^{(M)}$ has birth rates $b^{(M)}_1,\ldots,b^{(M)}_N>0$ and death rates $d^{(M)}_1,\ldots,d^{(M)}_M>0$, $d^{(M)}_{M+1}=\cdots=d^{(M)}_N=0$; see Figure~\ref{fig:twine} for a picture. In particular, we will choose $G^{(N-1)\,+}:=G$ and setting $G^+:=G^{(0)\,+}$ will yield the desired pure birth process with birth rates $b^+_1=\la_N,\ldots,b^+_N=\la_1$.

The core the proof is the following proposition, which corresponds to the inductive step in the argument.

\bp{\bf(Inductive step)}\label{P:indplus}
Let $1\leq M\leq N-1$ and let $G$ be the generator of a birth-and-death process in $\{0,\ldots,N\}$ with birth rates $b_1,\ldots,b_N>0$ and death rates $d_1,\ldots,d_M>0$, $d_{M+1}=\cdots=d_N=0$. Then there exists a probability kernel $K$ on $\{0,\ldots,N\}$ satisfying
\be
K(x,\{0,\ldots,x\})=1\quad(0\leq x\leq N)
\quad\mbox{and}\quad
K(x,x)=1\quad(M+1\leq x\leq N),
\ee
and a generator $G'$ of a birth-and-death process in $\{0,\ldots,N\}$ with birth rates $b'_1,\ldots,b'_N>0$ and death rates $d'_1,\ldots,d'_{M-1}>0$, $d'_M=\cdots=d'_N=0$, such that $KG=G'K$.
\ep
{\bf Proof} It follows from Lemma~\ref{L:eigen} applied to the process stopped at $M+1$ that there exists a function $\rho:\{0,\ldots,N\}\to\R$ such that $\rho>0$ on $\{0,\ldots,M\}$, $\rho=0$ on $\{M+1,\ldots,N\}$, $\sum_{x=0}^N\rho(x)=1$, and
\be\label{quasi}
G^\dgg\rho(x)=-\la\rho(x)+\la\de_{M+1}(x)\qquad(0\leq x\leq N),
\ee
where
\be
\la=b_{M+1}\rho(M)>0.
\ee
The law $\rho$ is sometimes called a {\em quasi-stationary law}. Using $\rho$, we define the kernel $K$ on $\{0,\ldots,N\}$ by
\be\label{Kplus1}
K(x,y):=\left\{\ba{ll}
\dis1_{\{y\leq x\}}\frac{\rho(y)}{H(x)}\quad&\mbox{if }x\leq M,\\[8pt]
\dis1_{\{y=x\}}\quad&\mbox{if }M+1\leq x,
\ea\right.
\ee
where
\be\label{Kplus2}
H(x):=\sum_{y=0}^x\rho(y)\qquad(0\leq x\leq M).
\ee
Since $K$ is a lower triangular matrix, it is invertible, so there exists a unique linear operator $G'$ satisfying $KG=G'K$ and $G'$ is in fact given by $G'=KGK^{-1}$. Since $G'1=G'K1=KG1=0$ we see that
\be
G'(x,x)=-\sum_{y\neq x}G'(x,y).
\ee
In view of this, to prove our claim, it suffices to check that the off-diagonal entries of $G'$ coincide with those of a birth-and-death process in $\{0,\ldots,N\}$ with birth rates $b'_1,\ldots,b'_N>0$ and death rates $d'_1,\ldots,d'_{M-1}>0$, $d'_M=\cdots=d'_N=0$.

To determine the off-diagonal entries of $G'$, we calculate, using (\ref{Gdagger}) and (\ref{quasi}),
\be\ba{l}\label{Kinv}
(KG)(x,y)=G^\dgg K(x,\,\cdot\,)(y)\\[5pt]
=\left\{\ba{ll}
\dis-\la1_{\{y\leq x\}}\frac{\rho(y)}{H(x)}-d_{x+1}\frac{\rho(x+1)}{H(x)}\de_x(y)+b_{x+1}\frac{\rho(x)}{H(x)}\de_{x+1}(y)\quad&\mbox{if }x<M,\\[5pt]
\dis-\la\rho(y)+\la\de_{M+1}(y)\quad&\mbox{if }x=M,\\[5pt]
\dis-b_{x+1}\de_x(y)+b_{x+1}\de_{x+1}(y)\quad&\mbox{if }x>M.
\ea\right.
\ec
In order to find $G'$, we need to express these formulas, as functions of $y$, as linear combinations of the basis vectors $(K(x,\,\cdot\,))_{0\leq x\leq N}$. To that aim, we observe that
\be
\de_x=K(x,\,\cdot\,)\qquad(M+1\leq x\leq N),
\ee
while for $1\leq x\leq M$, we have
\bc
\dis\de_x(y)
&=&\dis\big(1_{\{y\leq x\}}-1_{\{y\leq x-1\}}\big)\frac{\rho(y)}{\rho(x)}\\[5pt]
&=&\dis\frac{H(x)}{\rho(x)}1_{\{y\leq x\}}\frac{\rho(y)}{H(x)}-\frac{H(x-1)}{\rho(x)}1_{\{y\leq x-1\}}\frac{\rho(y)}{H(x-1)}\\[5pt]
&=&\dis\frac{H(x)}{\rho(x)}K(x,y)-\frac{H(x-1)}{\rho(x)}K(x-1,y).
\ec
Inserting this into (\ref{Kinv}), we find that
\be\ba{l}
\dis\sum_{x'=0}^NG'(x,x')K(x',y):=(KG)(x,y)\\[5pt]
=\left\{\ba{ll}
\dis-\la K(0,y)-d_1\frac{\rho(1)}{H(0)}K(0,y)\\[5pt]
\dis\quad+b_1\frac{\rho(0)}{H(0)}\Big(\frac{H(1)}{\rho(1)}K(1,y)-\frac{H(0)}{\rho(1)}K(0,y)\Big)\quad&\mbox{if }x=0,\\[5pt]
\dis-\la K(x,y)-d_{x+1}\frac{\rho(x+1)}{H(x)}\Big(\frac{H(x)}{\rho(x)}K(x,y)-\frac{H(x-1)}{\rho(x)}K(x-1,y)\Big)\\[5pt]
\dis\quad+b_{x+1}\frac{\rho(x)}{H(x)}\Big(\frac{H(x+1)}{\rho(x+1)}K(x+1,y)-\frac{H(x)}{\rho(x+1)}K(x,y)\Big)\quad&\mbox{if }0<x<M,\\[10pt]
\dis-\la K(M,y)+\la K(M+1,y)\quad&\mbox{if }x=M,\\[5pt]
\dis-b_{x+1}K(x,y)+b_{x+1}K(x+1,y)\quad&\mbox{if }x>M.
\ea\right.
\ec
From this, we can read off the off-diagonal entries of $G'$. Indeed,
\bc\label{bdac}
\dis b'_{x+1}=G'(x,x+1)&=&\left\{\ba{ll}
\dis b_{x+1}\frac{\rho(x)H(x+1)}{H(x)\rho(x+1)}\quad&\mbox{if }x<M,\\
\dis\la\quad&\mbox{if }x=M,\\[5pt]
\dis b_{x+1}\quad&\mbox{if }x>M,
\ea\right.\\[28pt]
\dis d'_x=G'(x,x-1)&=&\left\{\ba{ll}
%\dis 0\quad&\mbox{if }x=0,\\
\dis d_{x+1}\frac{\rho(x+1)H(x-1)}{H(x)\rho(x)}\quad&\mbox{if }0<x<M,\\
\dis 0\quad&\mbox{if }x\geq M,
\ea\right.
\ec
and all other off-diagonal entries are zero.\qed

%This remark dropped since not clear how to define $H(M+1)$ nor easy to see
%formula OK for $x=M$, nor easy to see d'_x positive
%\noi
%{\bf Remark 1} Using the `flux' relation
%\be
%b_x\rho(x-1)-d_x\rho(x)=\la H(x-1)\qquad(1\leq x\leq M),
%\ee
%we can rewrite (\ref{bdac}) as
%\bc
%\dis b'_x&=&\left\{\ba{ll}
%\dis d_x\frac{H(x)}{H(x-1)}+\la\frac{H(x)}{\rho(x)}\quad&\mbox{if }1\leq x\leq M+1,\\
%\dis b_{x+1}\quad&\mbox{otherwise},
%\ea\right.\\[28pt]
%\dis d'_x&=&\left\{\ba{ll}
%\dis b_{x+1}\frac{H(x-1)}{H(x)}-\la\frac{H(x-1)}{\rho(x)}\quad&\mbox{if }1\leq x\leq M-1,\\
%\dis 0\quad&\mbox{otherwise.}
%\ea\right.
%\ec
%This generalizes the formulas on \cite[page~564]{DM09} to the case $\la>0$ (taking into account that $b_{x+1}$ is called $b(x)$ in that paper).\med

\noi
{\bf Remark} The proof of Proposition~\ref{P:indplus} is straightforward except for the clever choice of $K$ in (\ref{Kplus1})--(\ref{Kplus2}). For some motivation of this choice and the way the authors arrived at it we refer to \cite{DM09}.\med

\noi
Using Proposition~\ref{P:indplus} we can construct a sequence of kernels $K^{(N-1)\,+},\ldots,K^{(1)\,+}$ and generators $G^{(N-1),+},\ldots,G^{(0)\,+}$ satisfying the intertwining relations (\ref{KGMplus}), such that $G^+:=G^{(0)\,+}$ is a pure birth process with birth rates $b^+_1,\ldots,b^+_N>0$, say. It is now easy to see that the composed kernel
\be
K^+:=K^{(1)\,+}\cdots K^{(N-1)\,+}
\ee
satisfies $K^+(x,\{0,\ldots,x\})=1$ ($0\leq x\leq N$), $K^+(N,N)=1$ and $K^+G=G^+K^+$. It is straightforward to check that the eigenvalues of $G'$ are $-b^+_1,\ldots,-b^+_N,0$. Since $G=(K^+)^{-1}G^+K^+$, the operators $G$ and $G^+$ have the same spectrum. 

We claim that $b^+_1>\cdots>b^+_N>0$. To see this, recall from the proofs of Lemma~\ref{L:eigen} and Proposition~\ref{P:indplus} that $-b^+_M$ is the Perron-Frobenius eigenvalue of the process with generator $G^{(M)\,+}$ stopped at $M+1$. It follows from the intertwining relation (\ref{KGMplus}) that $-b^+_{M-1}$ is also an eigenvalue of this process, corresponding to a different eigenvector, hence by the Perron-Frobenius theorem, $b_{M-1}>b_M$.

%Check new rates with DP at least in case la=0. warn of definition change.

%It seems all arguments go through if some of the death rates are zero. The only thing in Perron-Frobenius that may no longer be true is that all other eigenvalues are strictly smaller in value than $-\la$.

\subsection{Intertwining the slow process}

In the previous section, we have constructed a kernel $K^+$ and generator of a pure birth process $G^+$ such that (\ref{Gpm})~(i) holds. In this section, we construct a kernel $K^-$ and generator of a pure birth process $G^-$ satisfying (\ref{Gpm})~(ii). The proof will be very similar to the previous case, except that some things will `go he other way around'. In particular, using terminology introduced at the end of Section~\ref{S:twine}, $G^-$ will be the generator of an avaraged Markov process $X^-$ on $X$ while in the previous section we constructed a pure birth process $X^+$ such that $X$ is an averaged Markov process on $X^+$.

As in the previous section, the kernel $K^-$ will be constructed as the concatention of an inductively defined sequence of kernels $K^{(1)\,-},\ldots,K^{(N-1)\,-}$. Associated with these kernels is a sequence of generators $G^{(1),-},\ldots,G^{(N-1)\,-}$ of birth-and-death processes in $\{0,\ldots,N\}$ satisfying the intertwining relations
\be\label{KGMmin}
G^{(M-1)\,-}K^{(M)\,-}=K^{(M)\,-}G^{(M)\,-}\qquad(1\leq M\leq N-1),
\ee
where the process with generator $G^{(M)}$ has birth rates $b^{(M)}_1,\ldots,b^{(M)}_N>0$ and death rates $d^{(M)}_1=\cdots=d^{(M)}_M=0$, $d^{(M)}_{M+1},\ldots,d^{(M)}_{N-1}>0$, and $d^{(M)}_N=0$. We again refer to Figure~\ref{fig:twine} for an illustration.

The core of the argument is the following proposition.

\bp{\bf(Inductive step)}\label{P:indmin}
Let $0\leq M\leq N-2$ and let $G$ be the generator of a birth-and-death process in $\{0,\ldots,N\}$ with birth rates $b_1,\ldots,b_N>0$ and death rates $d_1=\cdots=d_M=0$, $d_{M+1},\ldots,d_{N-1}>0$, and $d_N=0$. Then there exists a probability kernel $K$ on $\{0,\ldots,N\}$ satisfying
\be
K(x,\{0,\ldots,x\})=1\quad(0\leq x\leq N)
\quad\mbox{and}\quad
K(x,x)=1\quad(x\not\in\{M,\ldots,N-1\}),
\ee
and a generator $\dot{G}$ of a birth-and-death process in $\{0,\ldots,N\}$ with birth rates $\dot{b}_1,\ldots,\dot{b}_N>0$ and death rates $\dot{d}_1=\cdots=\dot{d}_{M+1}=0$, $\dot{d}_{M+2},\ldots,\dot{d}_{N-1}>0$, and $\dot{d}_N=0$, such that $GK=K\dot{G}$.
\ep
{\bf Proof} It follows from Lemma~\ref{L:eigen} applied to the process restricted to $\{M,\ldots,N\}$ that there exists a function $f:\{0,\ldots,N\}\to\R$ such that $f=0$ on $\{0,\ldots,M-1\}$, $f$ is strictly decreasing on $\{M,\ldots,N\}$, $f(M)=1$, $f(N)=0$, and
\be
Gf(x)=-\la f(x)+b_M\de_{M-1}(x)\qquad(0\leq x\leq N),
\ee
where
\be\label{laid}
\la=b_{M+1}\big(1-f(M+1)\big)>0.
\ee
We set
\be
K(x,y):=1_{\{x=y\}}\qquad(y\not\in\{M,\ldots,N-1\}).
\ee
For $y=M,\ldots,N-1$, we claim that we can inductively define the kernel $K(x,y)$ and contants $C_y>0$ in such a way that
\be\label{indC}
\left.\ba{rl}
{\rm(i)}&\dis K(x,y):=C_y1_{\{y\leq x\}}f(x),\\[5pt]
{\rm(ii)}&\dis\sum_{y'=M}^yK(y,y')=1,
\ea\ \right\}\quad(M\leq y\leq N-1).
\ee
To see that this is all right, note that for $y=M$ (\ref{indC})~(i) and (ii) are satisfied by choosing $C_M:=1$, while for $M+1\leq y\leq N-1$ (\ref{indC})~(i) and (ii) imply that we must choose
\be
C_y:=\frac{1}{f(y)}\Big(1-\sum_{y'=M}^{y-1}K(y,y')\Big).
\ee
Since $f$ is strictly decreasing on $\{M,\ldots,N\}$, one has, by induction,
\be
\sum_{y'=M}^{y-1}K(y,y')=\sum_{y'=M}^{y-1}C_{y'}f(y)<\sum_{y'=M}^{y-1}C_{y'}f(y-1)=\sum_{y'=M}^{y-1}K(y-1,y')=1,
\ee
which shows that $C_y>0$. We now calculate
\be\ba{l}\label{GKform}
\dis(GK)(x,y)=GK(\,\cdot\,,y)(x)\\[5pt]
\quad=\left\{\ba{ll}
\dis b_y\de_{y-1}(x)-b_{y+1}\de_y(x)
\qquad&\mbox{if }0\leq y\leq M-1,\\[5pt]
\dis-\la f(x)+b_M\de_{M-1}(x)
\qquad&\mbox{if }y=M,\\[5pt]
\dis-\la C_y1_{\{y\leq x\}}f(x)+b_yC_y\de_{y-1}(x)\\[3pt]
\dis\quad-d_yC_y\big(f(y-1)-f(y)\big)\de_y(x)
\qquad&\mbox{if }M+1\leq y\leq N-1,\\[5pt]
b_N\de_{N-1}(x)\qquad&\mbox{if }y=N.
\ea\right.
\ec
By the same arguments as those in the previous section, there exists a unique linear operator $\dot{G}$ such that $GK=K\dot{G}$. In order to check that $\dot{G}$ is the generator of a birth-and-death process in $\{0,\ldots,N\}$ with birth rates $\dot{b}_1,\ldots,\dot{b}_N>0$ and death rates $\dot{d}_1=\cdots=\dot{d}_{M+1}=0$, $\dot{d}_{M+2},\ldots,\dot{d}_{N-1}>0$, and $\dot{d}_N=0$, it suffices to check that the off-diagonal entries $\dot{G}(x,y)$ have the desired form. In order to do this, we must express the formulas in (\ref{GKform}), as functions of $x$, as linear combinations of the basis vectors $(K(\,\cdot\,,y))_{0\leq y\leq N}$. We observe that
\be
\de_y(x)=K(\,\cdot\,,y)\qquad\big(y\not\in\{M,\ldots,N-1\}\big),
\ee
while for $M\leq y\leq N-2$, we have
\bc
\dis\de_y(x)
&=&\dis\big(1_{\{y\leq x\}}-1_{\{y+1\leq x\}}\big)\frac{f(x)}{f(y)}\\[5pt]
&=&\dis\frac{1}{f(y)C_y}C_y1_{\{y\leq x\}}f(x)-\frac{1}{f(y)C_{y+1}}C_{y+1}1_{\{y+1\leq x\}}f(x)\\[5pt]
&=&\dis\frac{1}{f(y)C_y}K(x,y)-\frac{1}{f(y)C_{y+1}}K(x,y+1),
\ec
and
\be
\de_{N-1}(x)=\frac{1}{f(N-1)C_{N-1}}K(x,N-1).
\ee
Inserting this into (\ref{GKform}), we obtain
\be\ba{l}\label{dotcalc}
\dis\sum_{y'}K(x,y')\dot{G}(y',y)=(GK)(x,y)\\[5pt]
=\left\{\ba{ll}
\dis b_yK(x,y-1)-b_{y+1}K(x,y)
\qquad&\mbox{if }0\leq y\leq M-1,\\[5pt]
\dis-\la K(x,M)+b_MK(x,M-1)
\qquad&\mbox{if }y=M,\\[5pt]
\dis-\la K(x,y)+b_yC_y\Big(\frac{K(x,y-1)}{f(y-1)C_{y-1}}-\frac{K(x,y)}{f(y+1)C_y}\Big)\\[3pt]
\dis\quad-d_yC_y\big(f(y-1)-f(y)\big)\Big(\frac{K(x,y)}{f(y)C_y}-\frac{K(x,y+1)}{f(y)C_{y+1}}\Big)
\quad&\mbox{if }M+1\leq y\leq N-2,\\[5pt]
\dis-\la K(x,y)+b_yC_y\Big(\frac{K(x,y-1)}{f(y-1)C_{y-1}}-\frac{K(x,y)}{f(y+1)C_y}\Big)\\[3pt]
\dis\quad-d_yC_y\big(f(y-1)-f(y)\big)\frac{K(x,y)}{f(y)C_y}
\quad&\mbox{if }y=N-1,\\[5pt]
\dis b_N\frac{K(x,N-1)}{f(N-1)C_{N-1}}\qquad&\mbox{if }y=N,
\ea\right.
\ec
where we use the convention that $b_0=0$ and hence $b_0K(x,-1)=0$, regardless of the (fictive) value of $K(x,-1)$. From (\ref{dotcalc}) we can read off the off-diagonal entries of $\dot{G}$. Indeed,
\bc
\dis\dot{b}_y=\dot{G}(y-1,y)&=&\left\{\ba{ll}
\dis b_y\quad&\mbox{if }1\leq y\leq M,\\
\dis b_y\frac{C_y}{f(y-1)C_{y-1}}\quad&\mbox{if }M+1\leq y\leq N,
\ea\right.\\[28pt]
\dis\dot{d}_{y+1}=\dot{G}(y+1,y)&=&\left\{\ba{ll}
\dis 0\quad&\mbox{if }y\not\in\{M+1,\ldots,N-2\},\\
\dis d_y\frac{C_y(f(y-1)-f(y))}{C_{y+1}f(y)}\quad&\mbox{if }M+1\leq y\leq N-2,
\ea\right.
\ec
and all other off-diagonal entries are zero. We note that in particular, by (\ref{laid}) and the definition of the $C_y$'s,
\be\label{bla}
\dot{b}_{M+1}=\frac{b_{M+1}C_{M+1}}{f(M)C_M}=b_{M+1}C_{M+1}=b_{M+1}\big(1-f(M+1))=\la.
\ee
\qed

\noi
{\bf Remark} As in the case of Proposition~\ref{P:indplus}, the proof of Proposition~\ref{P:indmin} is straightforward except for the choice of the kernel $K$. We have guessed formula (\ref{indC}) by analogy with formula (\ref{Kplus1}), which is due to \cite{DM09}.\med

\noi
With the help of Proposition~\ref{P:indmin}, we can inductively define kernels  $K^{(1)\,-},\ldots,K^{(N-1)\,-}$ and operators $G^{(1),-},\ldots,G^{(N-1)\,-}$. Setting $G^-:=G^{(N-1)\,-}$ and
\be
K^-=K^{(1)\,-}\cdots K^{(N-1)\,-}
\ee
now yields a generator of a pure birth process with birth rates $b^-_1,\ldots,b^-_N$ and a kernel $K^-$ with the properties described in (\ref{Kprop})--(\ref{Gpm}).

In the same way as in the previous section, we see that $0,-b^-_1,\ldots,-b^-_N$ are the eigenvalues of $G$. To see that $0<b_1<\cdots<b_N$ we observe from (\ref{bla}) that $-b^-_M$ is the Perron-Frobenius eigenvalue of the process with generator $G^{(M)\,-}$ restricted to $\{M,\ldots,N\}$. It follows from the intertwining relation (\ref{KGMmin}) that $-b^+_{M+1}$ is also an eigenvalue of this process, corresponding to a different eigenvector, hence by the Perron-Frobenius theorem, $b_M<b_{M+1}$.

\subsection{Proof of the main theorem}

{\bf Proof of Theorem~\ref{T:main}} The existence of generators $G^-,G^+$ and kernels $K^-,K^+$ satisfying (\ref{Kprop})--(\ref{Gpm}) has been proved in the previous sections. By Proposition~\ref{P:twine}, it follows that $X^+$ and $X$ can be coupled such that (\ref{Xcoup})~(i) holds. By applying Proposition~\ref{P:twine} to the kernel $L$ from $\{0,\ldots,N\}^2$ to $\{0,\ldots,N\}$ given by
\be
L\big((x,y),z\big):=K^-(y,z)\qquad(0\leq x,y,z\leq N),
\ee
we see that $(X^+,X)$ and $X^-$ can be coupled in such a way that both (\ref{Xcoup})~(i) and (ii) hold.\qed

\vspace{1cm}

\parbox[t]{14cm}{\small
Jan M.~Swart\\
Institute of Information Theory and Automation
of the ASCR (\' UTIA)\\
Pod vod\'arenskou v\v e\v z\' i 4,
18208 Praha 8,
Czech Republic\\
e-mail: swart@utia.cas.cz}

\end{document}